\documentclass[12pt]{article}

\usepackage[enableskew]{youngtab}
\usepackage[leqno]{amsmath}
\usepackage{amsfonts}
\usepackage{amssymb}
\usepackage{amsthm}
\usepackage{amssymb}
\usepackage[all]{xy}
\usepackage{graphicx}
\usepackage{bbm}

\newcommand{\fant}{\mbox{ }}

\theoremstyle{plain}
\newtheorem{theorem}{Theorem}[section]
\newtheorem{corollary}[theorem]{Corollary}

\theoremstyle{definition}
\newtheorem{definition}[theorem]{Definition}
\newtheorem{example}[theorem]{Example}

\newcommand{\la}{\lambda}

\title{Codes and shifted codes of partitions}

\author{J. T. Hird\\
\small Department of Mathematics\\[-0.8ex]
\small North Carolina State University, North Carolina, USA\\
\small \texttt{jthird@ncsu.edu}\\
\and
Naihuan Jing\\
\small Department of Mathematics\\[-0.8ex]
\small North Carolina State University, North Carolina, USA\\
\small \texttt{jing@math.ncsu.edu}\\
\and
Ernest Stitzinger\\
\small Department of Mathematics\\[-0.8ex]
\small North Carolina State University, North Carolina, USA\\
\small \texttt{stitz@math.ncsu.edu}\\
}

\date{October 19, 2010}

\begin{document}
\maketitle

\begin{abstract}
In a recent paper, Carrell and Goulden found a combinatorial
identity of the Bernstein operators that they then used to prove
Bernstein's Theorem. We show that this identity is a straightforward
consequence of the classical result.  We also show how a similar
approach using the codes of partitions can be generalized from Schur
functions to also include Schur $Q$-functions and derive the
combinatorial formulation for both cases.  We then apply them by
examining the Littlewood-Richardson and Pieri Rules.
\end{abstract}

\section{Introduction}\label{sec:intro}

Let $\Lambda=\mathbb C[p_1, p_2, p_3\ldots]=\oplus_{n=0}^{\infty}
\Lambda_n$ be the ring of symmetric functions, where $p_n$ is the
power sum symmetric function of degree $n$. As a graded vector
space, $\Lambda$ has several linear bases such as the power sum
symmetric functions $p_{\lambda}$ and Schur functions $s_{\lambda}$
\cite{macdonald, stanley} indexed by partitions. One way to
construct Schur functions is to realize them as images of Bernstein
operators $B_n$, whose generating function $B(t)$ is a variant of
vertex operator \cite{zelevinsky}:
$$B(t)=\sum_{n\in \mathbb{Z}}B_nt^n=exp\left(\sum_{k\geq1}
\frac{t^k}{k}p_k\right)exp\left(-\sum_{k\geq1}t^{-k}\frac{\partial}{\partial
p_k}\right) $$ which acts on the space $\Lambda$. In this
construction the Schur function $s_{\lambda}$ is easily given by
$s_{(\la_1, \ldots, \la_l)}=B_{\la_1}\cdots B_{\la_l}\cdot1$. The
operator $B_n$ is a graded linear transformation of degree $n$
defined via its action on the power sum function $p_{\mu}$. We can
also define the operator $B_n$ on the basis of Schur functions. It
turns out the action of $B_n$ on Schur functions has a close
relationship with Maya diagrams \cite{DJKMO}, one of the oldest
configurations of partitions. Recently Carrell and Goulden \cite{cg}
have formulated the action of $B(t)$ in terms of codes of
partitions, which are certain combinatorial description of Maya
diagrams. Similar combinatorial structures have also been used in
Okounkov's work on random matrices \cite{okounkov}.

Carrell and Goulden use codes of partitions to compute the action of
Bernstein operators on the Schur function $s_\lambda$. They then use
this formula to prove Bernstein's Theorem, that $s_{(\la_1, \ldots,
\la_l)}=B_{\la_1}\cdots B_{\la_l}\cdot1$ and
$B_mB_n=-B_{n-1}B_{m+1}$. Their proof is combinatorial and they also
show that their identity can be used in Pl\"ucker relations and KP
hierarchies.

In this paper we will show that Carrel and Goulden's formula can be
easily obtained from the classical results using algebraic
properties satisfied by vertex operators. We will also generalize
the combinatorial structures to the case of Schur $Q$-functions and
derive a similar but simpler combinatorial formulation for the
associated vertex operator. The Schur $Q$-functions are certain
distinguished linear bases in the subring of symmetric functions:
$$\Lambda^-=\mathbb C[p_1, p_3, p_5 \ldots]=\oplus_{n=0}^{\infty}
\Lambda_n^-.$$

These symmetric functions were defined by I. Schur in his seminal
work \cite{schur} on projective representations of the symmetric
group $S_n$ (see also \cite{macdonald}).
 As pointed out in
\cite{j1} Schur functions and Schur $Q$-functions are two examples of
the celebrated Boson-Fermion correspondence, in which they can be
roughly viewed as untwisted and twisted pictures of the Fock space
representations respectively, and the vertex operators for Schur and
Schur $Q$-functions come from two different realizations of affine Lie
algebras. Taking the advantage of this grand picture we can give a
unified approach to derive the action of vertex operators on Schur
and Schur $Q$-functions.

First we can compute the action of Bernstein operator by using the
commutation relations:
\begin{equation}\label{comm1}
B_mB_n=-B_{n-1}B_{m+1}.
\end{equation}
The combinatorial structure of
codes then follows easily from the algebraic structure.

When we tensor the ring $\Lambda$ by the group algebra of
one-dimensional lattice $\mathbb Z$, the commutation relations
(\ref{comm1}) can be improved into the exact anti-commutation
relations of the vertex operators $X(t)$:
\begin{equation}\label{comm2}
X_mX_n=-X_nX_m,
\end{equation} thus we obtain our second and even
simpler proof of Carrell-Goulden's formula. Using the same idea we
can generalize this to the twisted Fock space $\Lambda^-=\mathbb
C[p_1, p_3, \ldots]$  and again we use the similar antisymmetry of
the components of the vertex operator $Y(z)$ (see \cite{j1}) to
study the action of the Schur $Q$-functions.

 We also formulate the action of the twisted vertex operators in
 terms of {\it shifted codes}. In this way we have unified codes
 and shifted codes in the context of vertex operators and
 Boson-Fermion correspondence.  We also show
how these combinatorial objects can help us derive the
Littlewood-Richardson Rule and the Pieri Rules.

\section{Codes of Partitions}\label{sec:codes}

Let $\lambda = (\lambda_1, \lambda_2, \ldots, \lambda_l)$ be a
decreasing sequence of positive integers, $\lambda_1 \geq \lambda_2
\geq \cdots \geq \lambda_l > 0$.  We say that $\lambda$ is a
\textit{partition} of $n$, denoted $\lambda \vdash n$, if $\lambda_1
+ \lambda_2 + \cdots + \lambda_l = n$.  We also say that the
\textit{weight} of the partition $\lambda$ is $|\lambda|=n$, and the
\textit{length} of the partition is $l(\lambda) = l$.

The \textit{Young diagram} of a partition $\lambda$ is the left-justified arrangement of
boxes with $\lambda_i$ boxes in the $i^{\text{th}}$ row from the top.  Since the parts of
$\lambda$ are weakly decreasing, the number of boxes in each row will be less than or equal
to the number of boxes in each row above it.

Define the \textit{code} of a partition $\lambda$ to be the doubly infinite sequence of
letters R and U obtained from the Young diagram of shape $\lambda$ as follows.  Consider
the Young diagram top and left aligned in the $4^{\text{th}}$ quadrant of the $xy$-plane
together with the negative $y$-axis and the positive $x$-axis.  Trace up the negative $y$-axis
to the bottom of the Young diagram, then along the bottommost edge of the Young diagram,
then right along the positive $x$-axis.  The code of the partition is the sequence of R's
and U's obtained from this path, where R corresponds to a unit right step and U corresponds
to a unit up step.

\begin{example}
Let $\lambda = (4,2,1)$.  Then the Young diagram of shape $\lambda$ in the $4^{\text{th}}$
quadrant of the $xy$-plane is shown below, with the path described above in bold.
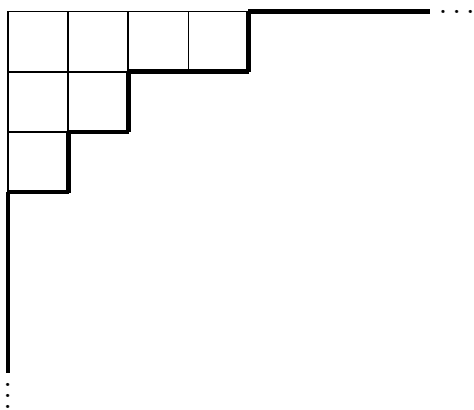
\begin{center}
\hspace{6.4cm}\dots\\
\vspace{-.05cm}
\setlength{\unitlength}{.8cm}
\begin{picture}(7,6)
\linethickness{0.075mm}
\put(0,6){\line(1,0){7}}
\put(0,5){\line(1,0){4}}
\put(0,4){\line(1,0){2}}
\put(0,3){\line(1,0){1}}
\put(0,6){\line(0,-1){6}}
\put(1,6){\line(0,-1){3}}
\put(2,6){\line(0,-1){2}}
\put(3,6){\line(0,-1){1}}
\put(4,6){\line(0,-1){1}}
\linethickness{.5mm}
\put(0,0){\line(0,1){3}}
\put(0,3){\line(1,0){1}}
\put(1,3){\line(0,1){1}}
\put(1,4){\line(1,0){1}}
\put(2,4){\line(0,1){1}}
\put(2,5){\line(1,0){2}}
\put(4,5){\line(0,1){1}}
\put(4,6){\line(1,0){3}}
\end{picture}\\
\vspace{-.1cm}
\hspace{-5.74cm}\fant\vdots
\end{center}
The path consists of infinitely many U's at the beginning - corresponding to tracing up the
negative $y$-axis, then RURURRU - corresponding to tracing the bottommost edge of the Young
diagram, then infinitely many R's at the end - corresponding to tracing right along the
positive $x$-axis.  Thus the sequence \ldots UUURURURRURRR \ldots is the code of the
partition $\lambda = (4,2,1)$.
\end{example}

Note that the code of any partition will always have infinitely many U's at the beginning
of the code, and infinitely many R's at the end of the code (corresponding respectively to
the negative $y$-axis and the positive $x$-axis).

Define the partition $\lambda^{(i)}$ to be the partition obtained by turning the
$i^{\text{th}}$ R from the left in the code of $\lambda$ to a U.  Equivalently, $\lambda^{(i)}$
is the partition obtained by looking at the lower-right edge of the associated Young diagram
(together with the positive $x$-axis and negative $y$-axis, where the Young diagram is
considered to be in the $4^{\text{th}}$ quadrant) and turning the $i^{\text{th}}$
horizontal edge from the left into a vertical edge (and shifting the resulting path
into the $4^{\text{th}}$ quadrant).

\begin{example}
Let $\lambda = (4,2,1)$.  Then to find $\lambda^{(3)}$, the third right step from the
left becomes an up step.  The changed edge is shown in bold below.

\vspace{.5cm}


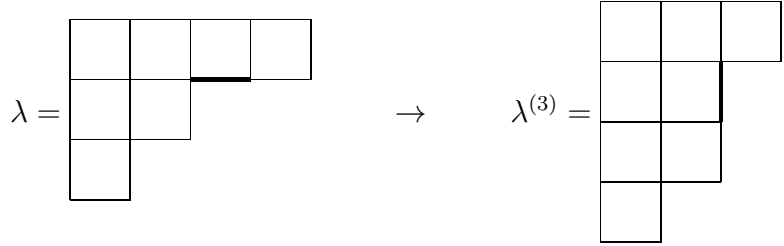
\begin{center}
{\Yvcentermath1 $\lambda =
\setlength{\unitlength}{.8cm}
\begin{picture}(4,2)
\linethickness{0.075mm}
\put(0,1.7){\line(1,0){4}}
\put(0,.7){\line(1,0){4}}
\put(0,-.3){\line(1,0){2}}
\put(0,-1.3){\line(1,0){1}}
\put(0,1.7){\line(0,-1){3}}
\put(1,1.7){\line(0,-1){3}}
\put(2,1.7){\line(0,-1){2}}
\put(3,1.7){\line(0,-1){1}}
\put(4,1.7){\line(0,-1){1}}
\linethickness{.5mm}
\put(2,.7){\line(1,0){1}}
\end{picture}
\hspace{1cm} \rightarrow \hspace{1cm} \lambda^{(3)} =
\setlength{\unitlength}{.8cm}
\begin{picture}(4,2)
\linethickness{0.075mm}
\put(0,2){\line(1,0){3}}
\put(0,1){\line(1,0){3}}
\put(0,0){\line(1,0){2}}
\put(0,-1){\line(1,0){2}}
\put(0,-2){\line(1,0){1}}
\put(0,2){\line(0,-1){4}}
\put(1,2){\line(0,-1){4}}
\put(2,2){\line(0,-1){3}}
\put(3,2){\line(0,-1){1}}
\linethickness{.5mm}
\put(2,0){\line(0,1){1}}
\end{picture}
$}
\end{center}

\vspace{1.8cm}
\noindent Thus $\lambda^{(3)} = (3,2,2,1)$.
\end{example}

A simple formula for the partition $\lambda^{(i)}$ is
\begin{eqnarray}\label{lambdi}
\lambda^{(i)} = (\lambda_1 - 1, \lambda_2 - 1, \ldots, \lambda_j - 1, i - 1, \lambda_{j+1}, \ldots, \lambda_l),
\end{eqnarray}
where $\lambda_j \geq i > \lambda_{j+1}$ (with the convention that $\lambda_{l+1}=0$ and
$\lambda_0=\infty$).

Given the code of a partition $\lambda$, let $u_i(\lambda)$ be the number of U's in the
code of $\lambda$ to the right of the $i^{\text{th}}$ R from the left, and let
 $r_i(\lambda)$ be the number of R's in the code of $\lambda$ to the left of the $i^{\text{th}}$ U
 from the right.  This means that $u_i(\lambda)$ is equal to the number of parts of $\lambda$ of
 size at least $i$, and $r_i(\lambda) = \lambda_i$.

\section{Bernstein Operators}\label{sec:bern}

Recall that the Bernstein operators $B(t),$ and $B_n$ are given by
$$B(t)=\sum_{n\in \mathbb{Z}}B_nt^n=exp\left(\sum_{k\geq1}\frac{t^k}{k}p_k\right)
exp\left(-\sum_{k\geq1}t^{-k}\frac{\partial}{\partial p_k}\right),
$$
which acts on the ring of polynomials $\Lambda=\mathbb C[p_1, p_2,
p_3, \ldots]$. Here the power sum $p_n$ acts as a multiplication on
$\Lambda$. Bernstein's primary result with these operators was
Bernstein's formula, which states:
\begin{eqnarray}\label{bern1}
s_\lambda = B_{\lambda_1}B_{\lambda_2}\cdots B_{\lambda_l}\cdot1,
\end{eqnarray}
where $\lambda = (\lambda_1, \lambda_2, \ldots, \lambda_l)$, and
$s_\lambda$ is the Schur polynomial indexed by $\lambda$.  For
convenience, we will often denote this composition as
$B_{\lambda_1}B_{\lambda_2}\cdots
B_{\lambda_l}\cdot1=B_{\lambda_1,\lambda_2,\ldots,\lambda_l}\cdot1$.
Another key relation satisfied by Bernstein operators is the
following:
\begin{eqnarray}\label{bern2}
B_nB_m=-B_{m-1}B_{n+1}.
\end{eqnarray}

We will now use these two results, Eqs.(\ref{bern1}) and
(\ref{bern2}), to prove a formula given in \cite{cg} which gives the
action of Bernstein's operators on the Schur polynomials.

\begin{theorem}\label{bern3}
For any partition $\lambda$,
$$B(t)s_\lambda=\sum_{i\geq1}(-1)^{|\lambda|-|\lambda^{(i)}|+i-1}t^{|\lambda^{(i)}|-|\lambda|}
s_{\lambda^{(i)}}.$$
\end{theorem}
This result was originally proved in \cite{cg} using some
combinatorial considerations and the dual action of the Schur
functions. We would like to give two simpler proofs to this result,
which will motivate our later generalization to the case of the Schur
$Q$-functions.

\textbf{Proof:} Since $B(t)=\sum_{n\in\mathbb{Z}}B_nt^n$, we only
need to determine the action of $B_n$ on $s_\lambda$.  By equation
(\ref{bern1}),
$$B_ns_\lambda=B_nB_{\lambda_1}B_{\lambda_2}\cdots B_{\lambda_l}\cdot1.$$

\textbf{Case 1}: If $n\geq\lambda_1$, then by equation (\ref{bern1}),
$$B_{n,\lambda_1,\lambda_2,\ldots,\lambda_l}\cdot1=s_{(n,\lambda_1,\lambda_2,\ldots,\lambda_l)}
=s_{\lambda^{(n+1)}},$$
where this term in the summation on the right has a $t$ term of $|\lambda^{(n+1)}|-|\lambda|$.
Since $n\geq\lambda_1$, turning the $(n+1)^\text{th}$ horizontal edge to a vertical edge creates
a new first row of size $n$.  So $|\lambda^{(n+1)}|-|\lambda|=n$ which is the exponent of $t$
associated with $B_n$.\\

\textbf{Case 2}: If $n=\lambda_j-j$ for some $j$, $1\leq j\leq l$, then by equation (\ref{bern2}),
\begin{eqnarray*}
B_{n,\lambda_1,\lambda_2,\ldots,\lambda_l}\cdot1 & = & (-1)\phantom{|}B_{\lambda_1-1,n+1,
\lambda_2,\ldots,\lambda_l}\cdot1\\
 & = & (-1)^2\phantom{|}B_{\lambda_1-1,\lambda_2-1,n+2,\lambda_3,\ldots,\lambda_l}\cdot1\\
 & \vdots &\\
 & = & (-1)^{j-1}\phantom{|}B_{\lambda_1-1,\lambda_2-1,\ldots,\lambda_{j-1}-1,n+j-1,\lambda_j,
 \lambda_{j+1},\ldots,\lambda_l}\cdot1,
\end{eqnarray*}
but $n=\lambda_j-j$, so $n+j-1=\lambda_j-1$.  From equation (\ref{bern2}), $B_{i,i+1}=-B_{i,i+1}$
for all $i$, which implies $B_{i,i+1}=0$ for all $i$.  Since $B_{n+j-1,\lambda_j}=B_{\lambda_j-1,
\lambda_j}$ is such a term, this product is zero.\\

\textbf{Case 3}: If $\lambda_{j+1}-(j+1)<n<\lambda_j-j$ for some $j$, $1\leq j<l$, then
similarly,
\begin{eqnarray*}
B_{n,\lambda_1,\lambda_2,\ldots,\lambda_l}\cdot1 & = & (-1)^{j}\phantom{|}B_{\lambda_1-1,
\lambda_2-1,\ldots,\lambda_j-1,n+j,\lambda_{j+1},\ldots,\lambda_l}\cdot1\\
 & = & (-1)^{j}\phantom{|}s_{(\lambda_1-1,\lambda_2-1,\ldots,\lambda_j-1,n+j,\lambda_{j+1},
 \ldots,\lambda_l)}\\
 & = & (-1)^{(-n)+(n+j+1)-1}\phantom{|}s_{\lambda^{(n+j+1)}}\\
 & = & (-1)^{(|\lambda|-|\lambda^{(n+j+1)}|)+(n+j+1)-1}\phantom{|}s_{\lambda^{(n+j+1)}},
\end{eqnarray*}
by equation (\ref{bern1}), since $\lambda_1-1\geq\lambda_2-1\geq\cdots\geq\lambda_j-1\geq
n+j\geq\lambda_{j+1}\geq\cdots\geq\lambda_l$.  Note that $|\lambda|-|\lambda^{(n+j+1)}|=-n$
since $\lambda^{(n+j+1)}$ removes the last box from each of the first $j$ rows of $\lambda$'s
Young diagram and then adds a row of size $n+j$.  Also note that the exponent of $t$
associated with $s_{\lambda^{(n+j+1)}}$ is $|\lambda^{(n+j+1)}|-|\lambda|=n$, the same
exponent associated with $B_n$.\\

\textbf{Case 4}: If $n<\lambda_l-l$, then similarly,
$$B_{n,\lambda_1,\lambda_2,\ldots,\lambda_l}\cdot1=(-1)^{l}\phantom{|}B_{\lambda_1-1,
\lambda_2-1,\ldots,\lambda_l-1,n+l}\cdot1.$$
$\bullet$ If $n+l\geq 0$, then by equation (\ref{bern1}),
\begin{eqnarray*}
(-1)^{l}\phantom{|}B_{\lambda_1-1,\lambda_2-1,\ldots,\lambda_l-1,n+l}\cdot1 & = &
(-1)^{(-n)+(n+l+1)-1}\phantom{|}s_{\lambda^{(n+l+1)}}\\
& = & (-1)^{|\lambda|-|\lambda^{(n+l+1)}|+(n+l+1)-1}\phantom{|}s_{\lambda^{(n+l+1)}},
\end{eqnarray*}
because $|\lambda|-|\lambda^{(n+l+1)}|=-n$ since $\lambda^{(n+l+1)}$ removes the last
box from each of the $l$ rows of $\lambda$'s Young diagram and then adds a row of size
$n+l$.  Again note that the exponent of $t$ associated with $s_{\lambda^{(n+l+1)}}$ is
$|\lambda^{(n+l+1)}|-|\lambda|=n$, the same exponent associated with $B_n$.\\
$\bullet$ If $n+l<0$, then by equation (\ref{bern2}):
$$B_{-1,0}=-B_{-1,0}=0$$
$$B_{-a,0}=-B_{-1,-a+1}=B_{-1,-1,-a+2}=\cdots= (-1)^a \phantom{|} B_{-1,-1,\ldots,-1,0}=0,$$
for all $a\in\mathbb{Z^+}$, since $B_0 \cdot 1 = 1$.  This implies that
$$(-1)^{l}\phantom{|}B_{\lambda_1-1,\lambda_2-1,\ldots,\lambda_l-1,n+l}\cdot1=(-1)^{l}
\phantom{|}B_{\lambda_1-1,\lambda_2-1,\ldots,\lambda_l-1,-a}\cdot1=0.$$
This proves the theorem.
\hfill $\Box$

We can also prove this theorem using vertex operators.  This method will be particularly
interesting to us because the same approach can be used to analyze the Schur $Q$-functions.

To see the symmetry of the indices of the Schur functions, we use a modified version of
Bernstein's operator from \cite{j1}.  Let $\mathbb{C}[\mathbb{Z}]$ be the group algebra
of $\mathbb{Z}$ generated by $e^p$, meaning $\mathbb{C}[\mathbb{Z}] = \oplus_{n\in\mathbb{Z}}
\hspace{.05cm} \mathbb{C}e^{np}$.  Consider the two operators $e^p$ and $t^{\partial_p}$ on
$\mathbb{C}[\mathbb{Z}]$ defined by
\begin{eqnarray*}
e^p\cdot e^{np} & = & e^{(n+1)p}\\
t^{\partial_p}\cdot e^{np} & = & t^n e^{np}.
\end{eqnarray*}

Following [J1], the vertex operator $X(t)$ is defined on $\Lambda\otimes \mathbb{C}[\mathbb{Z}]$ by
$$X(t) = B(t^{-1})e^p t^{\partial_p} = \sum_{n\in\mathbb{Z}} X_nt^{-n}.$$

The following result was proved in \cite{j2}: the product of the
vertex operator $X(t)$ is antisymmetric, so $X_nX_m=-X_mX_n$, and we
have the following theorem, which is a modified version from
\cite{j2}.

\begin{theorem}
\begin{enumerate}
\item For any $l\in\mathbb{N}$, one has
$$X_{t_1}\cdots X_{t_l} = (-1)^{l(\sigma)} X_{t_{\sigma(1)}}\cdots X_{t_{\sigma(l)}},$$
for all $\sigma$ in $S_{l}$, where $l(\sigma)$ is the number of inversions in the permutation
$\sigma$.\\
\item For any partition $\mu = (\mu_1, \ldots, \mu_l)$, we have
$$X_{-\mu_1} \cdots X_{-\mu_l} \cdot e^{mp} = s_{\mu-\delta+l\mathbbm{1}}e^{(m+l)p},$$
where $\delta = (l-1,\ldots, 2, 1, 0)$ and $\mathbbm{1}=(1, \ldots, 1) \in \mathbb{N}^l$.
\end{enumerate}
\end{theorem}
In particular, this means that
$$X_{-\mu_1} \cdots X_{-\mu_l} \cdot e^{-lp} = s_{\mu-\delta+l\mathbbm{1}}.$$
For simplicity, we removed the index shift of $\frac{1}{2}$ in the
definition of $X(t)$ (see \cite{j1, j2}).

Now we can give a simpler proof of Theorem \ref{bern3}.  For
simplicity, we will denote the composition as
$X_{-\mu_1}X_{-\mu_2}\cdots
X_{-\mu_l}=X_{-\mu_1,-\mu_2,\ldots,-\mu_l}$.
\begin{eqnarray*}
B(t)s_\lambda & = & B(t)X_{-(\lambda_1-1), -(\lambda_2-2), \ldots, -(\lambda_l-l)}\cdot e^{-lp}\\
& = & X(t)(e^pt^{\partial_p})^{-1}X_{-(\lambda_1-1), -(\lambda_2-2), \ldots, -(\lambda_l-l)}\cdot e^{-lp}\\
& = & X(t)X_{-(\lambda_1-2), -(\lambda_2-3), \ldots, -(\lambda_l-l-1)}\cdot (e^p t^{\partial_p})^{-1}t^{-l}e^{-lp}\\
& = & X(t)X_{-(\lambda_1-2), -(\lambda_2-3), \ldots, -(\lambda_l-l-1)}\cdot e^{-(l+1)p}\\
& = & \sum_{n\in \mathbb{Z}}X_{-n}X_{-(\lambda_1-2), -(\lambda_2-3), \ldots,
-(\lambda_l-l-1)}\cdot e^{-(l+1)p}t^n\\
& = & \sum_{n\neq \lambda_k-k-1}(-1)^jX_{-(\lambda_1-2), \ldots, -(\lambda_j-j-1), -n,
-(\lambda_{j+1}-j-2), \ldots, -(\lambda_l-l-1)}\cdot e^{-(l+1)p}t^n\\
& = & \sum_{n\neq \lambda_k-k-1}(-1)^jX_{-(\lambda^{(i)}_1-1), \ldots,
-(\lambda^{(i)}_{l+1}-l-1)}\cdot e^{-(l+1)p}t^n
s_{\lambda^{(i)}},
\end{eqnarray*}
where $\lambda_j-j-1 > n > \lambda_{j+1}-j-2$ and $i=n+j+1$, so $\lambda_j>i\geq\lambda_{i+1}-1$.
This definition of $i$ also implies that $n=|\lambda^{(i)}|-|\lambda|$,
$j=|\lambda|-|\lambda^{(i)}|+i-1$, and $\lambda^{(i)}=(\lambda_1-1, \ldots, \lambda_{j}-1, i-1,
\lambda_{j+1}, \ldots, \lambda_l) = (\lambda_1-1, \ldots, \lambda_{j}-1, n+j, \lambda_{j+1},
\ldots, \lambda_l)$.  With this identification this last line becomes the following:
$$B(t)s_\lambda=\sum_{i\geq1}(-1)^{|\lambda|-|\lambda^{(i)}|+i-1}t^{|\lambda^{(i)}|-
|\lambda|}s_{\lambda^{(i)}}$$
and Theorem \ref{bern3} is proved.
\hfill $\Box$

\section{Schur $Q$-function}\label{sec:schQ}

We will next state and prove a similar result for the Schur
$Q$-functions, $Q_{\lambda}$, where $\lambda$ is a \textit{strict
partition}, i.e. $\lambda_1>\lambda_2>\cdots>\lambda_l$ and
$\la_i\in\mathbb N$.

For any partition $\mu=(1^{m_1(\mu)}2^{m_1(\mu)}\cdots)$, we define
$z_{\mu}=\prod_{i\geq 1}i^{m_i(\mu)}m_i(\mu)!$. We consider the ring
of symmetric functions in $x_1, x_2, \ldots$, but restrict ourselves
to polynomials in odd degree power sums
$$
p_{2k+1}=\sum_{i\geq 1}x_i^{2k+1}, \qquad k\in\mathbb Z_+.
$$
Let $\mathcal OP$ denote the set of partitions with odd parts, and
let $\Lambda^-$ be the ring of symmetric functions generated by
$p_{2k+1}, k\in\mathbb Z_+$. Under the inner product
$$<p_{\la},
p_{\mu}>=2^{-l(\la)}\delta_{\la, \mu}z_{\la}, \qquad
\la,\mu\in\mathcal{OP}, $$ the space $\Lambda^-$ has $Q_{\la}$
($\la$ strict) as a distinguished orthogonal basis of symmetric
polynomials \cite{schur, macdonald}. They play a fundamental role in
the construction of projective representations of the symmetric
group $S_n$.

On the space $\Lambda^-$ we recall the definition of the twisted
vertex operator \cite{j1}:
$$Y(t)=\sum_{n\in \mathbb{Z}}Y_nt^{-n}=exp\left(\sum_{k\geq1}\frac{2t^{-2k+1}}{k}p_{2k-1}\right)
exp\left(-\sum_{k\geq1}t^{2k-1}\frac{\partial}{\partial
p_{2k-1}}\right),
$$
which acts on the ring of polynomials $\Lambda^-=\mathbb C[p_1, p_3,
p_5, \ldots]$, and the power sum $p_{2k-1}$ acts as a multiplication
on $\Lambda^-$.

From \cite{j1}, we have that the following two results hold:
\begin{eqnarray}\label{jing1}
Q_\lambda = Y_{-\lambda_1}Y_{-\lambda_2}\cdots Y_{-\lambda_l}\cdot1,
\end{eqnarray}
where $\lambda = (\lambda_1, \lambda_2, \ldots, \lambda_l)$, and $Q_\lambda$ is the Schur
$Q$-function indexed by $\lambda$.  Again, we will often denote this composition as
$Y_{\lambda_1}Y_{\lambda_2}\cdots Y_{\lambda_l}\cdot1=Y_{\lambda_1,\lambda_2,\ldots,\lambda_l}
\cdot1$.  The second result is:
\begin{eqnarray}\label{jing2}
Y_nY_m=-Y_{m}Y_{n}.
\end{eqnarray}

\begin{theorem}\label{jing3}
For any strict partition $\lambda$,
$$Y(t)Q_\lambda=\sum_{n\neq\lambda_j}(-1)^{i}t^{n}Q_{(\lambda_1,\lambda_2,\ldots,\lambda_i,n,
\lambda_{i+1},\ldots,\lambda_l)}.$$
\end{theorem}

\textbf{Proof:} Recall that $Y(t)=\sum_{n\in\mathbb{Z}}Y_nt^{-n}$,
then use equations (\ref{jing1}) and (\ref{jing2}):
\begin{eqnarray*}
Y(t)Q_\lambda=\sum_{n\in\mathbb{Z}}Y_{n}t^{-n}Q_\lambda & = &
\sum_{n\in\mathbb{Z}}t^nY_{-n}Y_{-\lambda_1}
Y_{-\lambda_2}\cdots Y_{-\lambda_l}\cdot1\\
& = &
\sum_{n\neq\lambda_j}(-1)^{i}t^nY_{-\lambda_1,-\lambda_2,-\ldots,-\lambda_i,-n,-\lambda_{i+1},
\ldots,-\lambda_l}\cdot1\\
& = & \sum_{n\neq\lambda_j}(-1)^{i}t^{n}Q_{(\lambda_1,\lambda_2,\ldots,\lambda_i,n,\lambda_{i+1},
\ldots,\lambda_l)},
\end{eqnarray*}
where $\lambda_i>n>\lambda_{i+1}$, because by equation
(\ref{jing2}), $Y_{-n}Y_{-\lambda_j}=Y_{-n}Y_{-n}=0$ if
$n=\lambda_j$. \hfill $\Box$
\\

We can also interpret the result in terms of codes of strict partitions, but we
first need to reinterpret how codes behave for strict partitions.

\begin{definition}
Define the partition $\lambda^{[i]}$ to be the partition obtained from the code
of a strict partition $\lambda$ by inserting a U between the $i^{\text{th}}$ pair
of consecutive R's (with the convention that three consecutive R's counts as two
pairs, four consecutive R's counts as three pairs, and so on).  Equivalently,
$\lambda^{[i]}$ is the partition obtained from the code of $\lambda$ by inserting
a U after the $i^{\textit{th}}$ R which is immediately followed by an R.
\end{definition}

\begin{example} For example, if $\lambda = (6,4,3,1)$, the first pair of consecutive
R's in the code of $\lambda$ is shown below in bold, with the new edge inserted between
them to get $\lambda^{[1]}$ also shown in bold.

\vspace{.2cm}

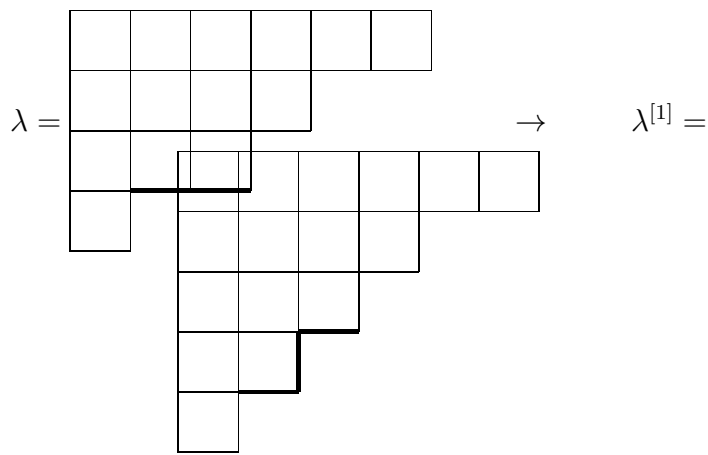
\begin{center}
{\Yvcentermath1 $\lambda =
\setlength{\unitlength}{.8cm}
\begin{picture}(6,3)
\linethickness{0.075mm}
\put(0,2){\line(1,0){6}}
\put(0,1){\line(1,0){6}}
\put(0,0){\line(1,0){4}}
\put(0,-1){\line(1,0){3}}
\put(0,-2){\line(1,0){1}}
\put(0,2){\line(0,-1){4}}
\put(1,2){\line(0,-1){4}}
\put(2,2){\line(0,-1){3}}
\put(3,2){\line(0,-1){3}}
\put(4,2){\line(0,-1){2}}
\put(5,2){\line(0,-1){1}}
\put(6,2){\line(0,-1){1}}
\linethickness{.5mm}
\put(1,-1){\line(1,0){2}}
\end{picture}
\hspace{1cm} \rightarrow \hspace{1cm} \lambda^{[1]} =
\setlength{\unitlength}{.8cm}
\begin{picture}(6,3)
\linethickness{0.075mm}
\put(0,2.7){\line(1,0){6}}
\put(0,1.7){\line(1,0){6}}
\put(0,.7){\line(1,0){4}}
\put(0,-.3){\line(1,0){3}}
\put(0,-1.3){\line(1,0){2}}
\put(0,-2.3){\line(1,0){1}}
\put(0,2.7){\line(0,-1){5}}
\put(1,2.7){\line(0,-1){5}}
\put(2,2.7){\line(0,-1){4}}
\put(3,2.7){\line(0,-1){3}}
\put(4,2.7){\line(0,-1){2}}
\put(5,2.7){\line(0,-1){1}}
\put(6,2.7){\line(0,-1){1}}
\linethickness{.5mm}
\put(1,-1.3){\line(1,0){1}}
\put(2,-.3){\line(1,0){1}}
\put(2,-1.3){\line(0,1){1}}
\end{picture}
$}
\end{center}

\vspace{2.4cm}

\noindent To get $\lambda^{[2]}$, we insert a U between the second pair of consecutive
R's in the code of $\lambda$.  Again the pair of right steps corresponding to those R's
are shown below in bold, along with the up step inserted between them.

\vspace{.2cm}

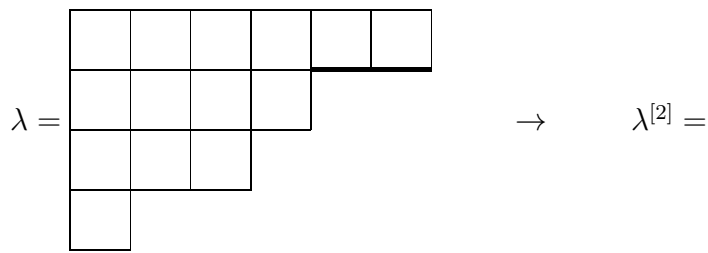
\begin{center}
{\Yvcentermath1 $\lambda =
\setlength{\unitlength}{.8cm}
\begin{picture}(6,3)
\linethickness{0.075mm}
\put(0,2){\line(1,0){6}}
\put(0,1){\line(1,0){6}}
\put(0,0){\line(1,0){4}}
\put(0,-1){\line(1,0){3}}
\put(0,-2){\line(1,0){1}}
\put(0,2){\line(0,-1){4}}
\put(1,2){\line(0,-1){4}}
\put(2,2){\line(0,-1){3}}
\put(3,2){\line(0,-1){3}}
\put(4,2){\line(0,-1){2}}
\put(5,2){\line(0,-1){1}}
\put(6,2){\line(0,-1){1}}
\linethickness{.5mm}
\put(4,1){\line(1,0){2}}
\end{picture}
\hspace{1cm} \rightarrow \hspace{1cm} \lambda^{[2]} =
\setlength{\unitlength}{.8cm}
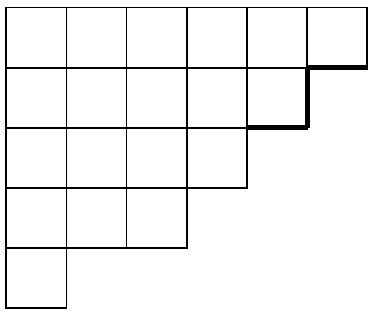
\begin{picture}(6,3)
\linethickness{0.075mm}
\put(0,2.7){\line(1,0){6}}
\put(0,1.7){\line(1,0){6}}
\put(0,.7){\line(1,0){5}}
\put(0,-.3){\line(1,0){4}}
\put(0,-1.3){\line(1,0){3}}
\put(0,-2.3){\line(1,0){1}}
\put(0,2.7){\line(0,-1){5}}
\put(1,2.7){\line(0,-1){5}}
\put(2,2.7){\line(0,-1){4}}
\put(3,2.7){\line(0,-1){4}}
\put(4,2.7){\line(0,-1){3}}
\put(5,2.7){\line(0,-1){2}}
\put(6,2.7){\line(0,-1){1}}
\linethickness{.5mm}
\put(4,.7){\line(1,0){1}}
\put(5,1.7){\line(1,0){1}}
\put(5,.7){\line(0,1){1}}
\end{picture}
$}
\end{center}

\vspace{2.4cm}

\noindent To get $\lambda^{[3]}$, we insert a U between the third pair of consecutive R's
in the code of $\lambda$.  This works the same way as the previous examples, except that the
third pair of R's are in the part of the code corresponding to the positive $x$-axis.  Again
the pair of right steps corresponding to those R's are shown below in bold, along with the up
step inserted between them.

\vspace{.2cm}

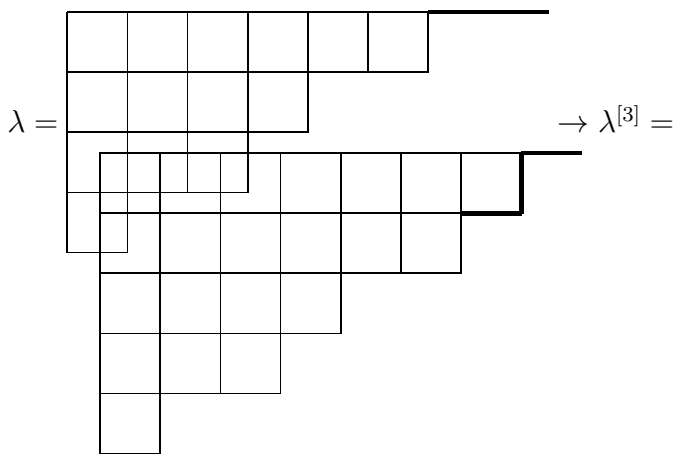
\begin{center}
{\Yvcentermath1 $\lambda =
\setlength{\unitlength}{.8cm}
\begin{picture}(8,3)
\linethickness{0.075mm}
\put(0,2){\line(1,0){8}}
\put(0,1){\line(1,0){6}}
\put(0,0){\line(1,0){4}}
\put(0,-1){\line(1,0){3}}
\put(0,-2){\line(1,0){1}}
\put(0,2){\line(0,-1){4}}
\put(1,2){\line(0,-1){4}}
\put(2,2){\line(0,-1){3}}
\put(3,2){\line(0,-1){3}}
\put(4,2){\line(0,-1){2}}
\put(5,2){\line(0,-1){1}}
\put(6,2){\line(0,-1){1}}
\linethickness{.5mm}
\put(6,2){\line(1,0){2}}
\end{picture}
\rightarrow \lambda^{[3]} =
\setlength{\unitlength}{.8cm}
\begin{picture}(8,3)
\linethickness{0.075mm}
\put(0,2.7){\line(1,0){8}}
\put(0,1.7){\line(1,0){7}}
\put(0,.7){\line(1,0){6}}
\put(0,-.3){\line(1,0){4}}
\put(0,-1.3){\line(1,0){3}}
\put(0,-2.3){\line(1,0){1}}
\put(0,2.7){\line(0,-1){5}}
\put(1,2.7){\line(0,-1){5}}
\put(2,2.7){\line(0,-1){4}}
\put(3,2.7){\line(0,-1){4}}
\put(4,2.7){\line(0,-1){3}}
\put(5,2.7){\line(0,-1){2}}
\put(6,2.7){\line(0,-1){2}}
\put(7,2.7){\line(0,-1){1}}
\linethickness{.5mm}
\put(6,1.7){\line(1,0){1}}
\put(7,2.7){\line(1,0){1}}
\put(7,1.7){\line(0,1){1}}
\end{picture}
$}
\end{center}

\end{example}

\vspace{2.4cm}

Another way to think about $\lambda^{[i]}$ is the following.  With this definition
$\lambda^{[i]}$ is the strict partition with the $i^{\text{th}}$ smallest possible integer
inserted into the partition $\lambda$.  This means that $\lambda^{[1]}$ is the strict partition
with the smallest possible integer inserted into $\lambda$.  For $\lambda = (6, 4, 3, 1),$
the smallest integer that can be inserted to still have a strict partition is 2,
so $\lambda^{[1]} = (6, 4, 3, 2, 1)$.  The second smallest integer that can be inserted
into $\lambda = (6, 4, 3, 1)$ is 5, so $\lambda^{[1]} = (6, 5, 4, 3, 1)$.  Similarly,
$\lambda^{[3]} = (7, 6, 4, 3, 1)$, $\lambda^{[4]} = (8, 6, 4, 3, 1)$, and so on.


Often strict partitions are associated with \textit{shifted Young
diagram} \cite{stembridge} rather than Young diagram.  A shifted
Young diagram of shape $\lambda$, where $\lambda$ is a strict
partition, is an arrangement of boxes with $\lambda_i$ boxes in the
$i^{\text{th}}$ row, with the leftmost box in each row one unit to
the right of the leftmost box of the row above it.  This is
sometimes more intuitive since the rightmost edge of a shifted Young
diagram of shape $\lambda$, where $\lambda$ is a strict partition,
follows the same rules of a Young diagram of shape $\mu$, where
$\mu$ is any partition, namely that the rightmost edge moves weakly
left as you go from top to bottom. We can use this correlation to
reinterpret $\lambda^{[i]}$ using the analogue of our existing
machinery for codes on a shifted Young diagram of shape $\lambda$.

\begin{definition}
Define the \textit{shifted code} of a strict partition $\lambda$ to
be the infinite sequence of letters R and U obtained from the
shifted Young diagram of shape $\lambda$ as follows.  Consider the
shifted Young diagram top and left aligned in the $4^{\text{th}}$
quadrant of the $xy$-plane together with the positive $x$-axis.
Starting at the bottom right corner of the leftmost box in the last
row, trace along the rightmost edge of the shifted Young diagram,
then right along the positive $x$-axis.  Equivalently, start the
code at the lowest place where the line $y=-x$ intersects the
shifted Young diagram.  The shifted code of the strict partition is
the sequence of R's and U's obtained from this path, where R
corresponds to a unit right step and U corresponds to a unit up
step.
\end{definition}

\begin{example}
Let $\lambda = (5,4,2)$.  Then the shifted Young diagram of shape
$\lambda$ in the $4^{\text{th}}$ quadrant of the $xy$-plane is shown
below, with the path described above in bold.
\begin{center}
\hspace{7.2cm}\dots\\
\vspace{-.05cm}
\setlength{\unitlength}{.8cm}
\begin{picture}(8,3)
\linethickness{0.075mm}
\put(0,3){\line(1,0){8}}
\put(0,2){\line(1,0){5}}
\put(1,1){\line(1,0){4}}
\put(2,0){\line(1,0){2}}
\put(0,3){\line(0,-1){1}}
\put(1,3){\line(0,-1){2}}
\put(2,3){\line(0,-1){3}}
\put(3,3){\line(0,-1){3}}
\put(4,3){\line(0,-1){3}}
\put(5,3){\line(0,-1){2}}
\linethickness{.5mm}
\put(3,0){\line(1,0){1}}
\put(4,0){\line(0,1){1}}
\put(4,1){\line(1,0){1}}
\put(5,1){\line(0,1){2}}
\put(5,3){\line(1,0){3}}
\end{picture}\\
\end{center}
\end{example}

\begin{example}
Let $\lambda = (6,4,3,1)$.  Then the shifted Young diagram of shape
$\lambda$ in the $4^{\text{th}}$ quadrant of the $xy$-plane is shown
below, with the path described above in bold.
\begin{center}
\hspace{8cm}\dots\\
\vspace{-.05cm}
\setlength{\unitlength}{.8cm}
\begin{picture}(9,4)
\linethickness{0.075mm}
\put(0,4){\line(1,0){9}}
\put(0,3){\line(1,0){6}}
\put(1,2){\line(1,0){4}}
\put(2,1){\line(1,0){3}}
\put(3,0){\line(1,0){1}}
\put(0,4){\line(0,-1){1}}
\put(1,4){\line(0,-1){2}}
\put(2,4){\line(0,-1){3}}
\put(3,4){\line(0,-1){4}}
\put(4,4){\line(0,-1){4}}
\put(5,4){\line(0,-1){3}}
\put(6,4){\line(0,-1){1}}
\linethickness{.5mm}
\put(4,0){\line(0,1){1}}
\put(4,1){\line(1,0){1}}
\put(5,1){\line(0,1){2}}
\put(5,3){\line(1,0){1}}
\put(6,3){\line(0,1){1}}
\put(6,4){\line(1,0){3}}
\end{picture}\\
\end{center}
\end{example}

Note that the shifted code is not doubly infinite like the code of
an arbitrary partition, since it has a fixed starting point.  It
does however retain the property that there are infinitely many R's
at the end of the code.

Using shifted codes we can reinterpret our definition of
$\lambda^{[i]}$.  For a strict partition $\lambda$, $\lambda^{[i]}$
is obtained from the shifted code of $\lambda$ by turning the
$i^{\text{th}}$ R in the shifted code to a U.  This is since either
method inserts the $i^{\text{th}}$ smallest possible integer into
the partition $\lambda$ to still have a strict partition, or since
the number of pairs of consecutive R's between two U's is the number
of consecutive R's minus one, which is the number of R's in the
shifted code corresponding to the same row.

\begin{example}
We return to our example $\lambda = (6,4,3,1)$.  Then we can find
$\lambda^{[2]}$ by turning the second right step from the left in
the shifted code of $\lambda$ into an up step.  The changed edge is
shown in bold below.

\vspace{.2cm}

\begin{center}
{\Yvcentermath1 $\lambda =
\setlength{\unitlength}{.8cm}
\begin{picture}(6,3)
\linethickness{0.075mm}
\put(0,2){\line(1,0){6}}
\put(0,1){\line(1,0){6}}
\put(1,0){\line(1,0){4}}
\put(2,-1){\line(1,0){3}}
\put(3,-2){\line(1,0){1}}
\put(0,2){\line(0,-1){1}}
\put(1,2){\line(0,-1){2}}
\put(2,2){\line(0,-1){3}}
\put(3,2){\line(0,-1){4}}
\put(4,2){\line(0,-1){4}}
\put(5,2){\line(0,-1){3}}
\put(6,2){\line(0,-1){1}}
\linethickness{.5mm}
\put(5,1){\line(1,0){1}}
\end{picture}
\hspace{1cm} \rightarrow \hspace{1cm} \lambda^{[2]} =
\setlength{\unitlength}{.8cm}
\begin{picture}(6,3)
\linethickness{0.075mm}
\put(0,2.7){\line(1,0){6}}
\put(0,1.7){\line(1,0){6}}
\put(1,.7){\line(1,0){5}}
\put(2,-.3){\line(1,0){4}}
\put(3,-1.3){\line(1,0){3}}
\put(4,-2.3){\line(1,0){1}}
\put(0,2.7){\line(0,-1){1}}
\put(1,2.7){\line(0,-1){2}}
\put(2,2.7){\line(0,-1){3}}
\put(3,2.7){\line(0,-1){4}}
\put(4,2.7){\line(0,-1){5}}
\put(5,2.7){\line(0,-1){5}}
\put(6,2.7){\line(0,-1){4}}
\linethickness{.5mm}
\put(6,.7){\line(0,1){1}}
\end{picture}
$}
\end{center}

\vspace{2.4cm}

\end{example}

Given the code of a strict partition $\lambda$, let $\tilde
u_i(\lambda)$ be the number of U's in the code of $\lambda$ to the
right of the $i^{\text{th}}$ pair of consecutive R's from the left,
which is the number of U's in the shifted code of $\lambda$ to the
right of the $i^{\text{th}}$ R from the left.

This means that $\tilde u_i(\lambda)$ is equal to the number of parts of $\lambda$
greater than the $i^{\text{th}}$ smallest possible integer that can be inserted into
$\lambda$, which is equal to the number of parts of $\lambda$ of size at least $|\lambda^{[i]}| - |\lambda|$.  Then the number of parts  of $\lambda$ of size at least $|\lambda^{[i]}| - |\lambda|$ is the length of $\lambda$ minus the number of parts of size less than $|\lambda^{[i]}| - |\lambda|$.  But the number of parts less than $|\lambda^{[i]}| - |\lambda|$ is the number of integers less than $|\lambda^{[i]}| - |\lambda|$ minus the number of integers less than $|\lambda^{[i]}| - |\lambda|$ that are not in $\lambda$, which is $(|\lambda^{[i]}| - |\lambda| - 1) - ( i - 1)$ = $|\lambda^{[i]}| - |\lambda| - i$.  So $\tilde u_i(\lambda) = l(\lambda) - (|\lambda^{[i]}| - |\lambda| - i) = l(\lambda) + |\lambda| - |\lambda^{[i]}| + i$.

We can now use $\lambda^{[i]}$ to reinterpret Theorem \ref{jing3}.

\begin{theorem}
For any strict partition $\lambda$,
$$Y(t)Q_\lambda=\sum_{i\geq1}(-1)^{l(\lambda) + |\lambda| - |\lambda^{[i]}|
+ i}t^{|\lambda^{[i]}| - |\lambda|} Q_{\lambda^{[i]}}.$$
\end{theorem}

\textbf{Proof:}  By Theorem \ref{jing3}, we know that
$$Y(t)Q_\lambda=\sum_{n\neq\lambda_j}(-1)^{k}t^{n}Q_{(\lambda_1,\lambda_2,\ldots,
\lambda_k,n,\lambda_{k+1},\ldots,\lambda_l)}.$$
But $\lambda^{[i]}$ is the partition with the $i^{\text{th}}$ smallest possible
integer that can be inserted into the partition $\lambda$.
Thus $(\lambda_1,\lambda_2,\ldots,\lambda_k,n,\lambda_{k+1},\ldots,\lambda_l)$ = $\lambda^{[i]}$
for some $i$, where $n$ is the $i^{\text{th}}$ smallest possible integer
that can be inserted into $\lambda$, so $n =  |\lambda^{[i]}| - |\lambda|$.
Then $k$ is the number of parts of $\lambda$ greater than the $i^{\text{th}}$
smallest possible integer that can be inserted into $\lambda$, so by definition
$k = \tilde u_i(\lambda)$.  Thus
\begin{eqnarray*}
Y(t)Q_\lambda & = & \sum_{n\neq\lambda_j}(-1)^{k}t^{n}
Q_{(\lambda_1,\lambda_2,\ldots,\lambda_k,n,\lambda_{k+1},\ldots,\lambda_l)}\\
 & = & \sum_{i\geq1}(-1)^{\tilde u_i(\lambda)}t^{|\lambda^{[i]}| - |\lambda|}
 Q_{\lambda^{[i]}}\\
 & = & \sum_{i\geq1}(-1)^{l(\lambda) + |\lambda| - |\lambda^{[i]}| + i}
 t^{|\lambda^{[i]}| - |\lambda|} Q_{\lambda^{[i]}}
\end{eqnarray*}
since we know $\tilde u_i(\lambda) = k = l(\lambda) + |\lambda| - |\lambda^{[i]}|
+ i$, $|\lambda^{[i]}| - |\lambda| = n$, and\\
$\lambda^{[i]} = (\lambda_1,\lambda_2,\ldots,\lambda_i,n,\lambda_{i+1},\ldots,\lambda_l)$.
\hfill $\Box$

\section{Littlewood-Richardson Rule}\label{sec:L-R}

One application for codes of partitions is the following theorem,
which gives a new way to compute Littlewood-Richardson coefficients
\cite{macdonald, sagan}, using only the codes of the partitions
involved.

A skew-partition $\lambda/\mu$ is a \textit{horizontal $n$-strip} if no column in
the Young diagram of $\lambda/\mu$ has more than one box.  Equivalently, $\lambda/\mu$
is a horizontal $n$-strip if $\lambda_{i+1} \leq \mu_i \leq \lambda_i$ for all
$1 \leq i \leq l(\mu)$, where $l(\mu)$ is the length of $\mu$.

A skew-partition $\lambda/\mu$ is a \textit{vertical $n$-strip} if no row in the
Young diagram of $\lambda/\mu$ has more than one box.  Equivalently, $\lambda/\mu$
is a horizontal $n$-strip if $\lambda_i - 1 \leq \mu_i \leq \lambda_i$ for
all $1 \leq i \leq l(\mu)$, where $l(\mu)$ is the length of $\mu$.

\begin{theorem}\label{l-r}
(The Littlewood-Richardson Rule)
$$s_\mu s_\nu = \sum_\lambda c^\lambda_{\mu,\nu} s_\lambda =
\sum_{(\mu=\mu^0, \mu^1, \mu^2, \ldots, \mu^l)} s_{\mu^l},$$
where $l = u_1(\nu)$.  Given the code of the partition $\mu^{i-1}$, $\mu^i$
is obtained as follows:
\begin{itemize}
\item Starting with the U left of the leftmost R in the code of $\lambda$ and
working to the right, move the U's to the right a total of $r_i(\nu)$ places
by switching a UR to RU in the code $r_i(\nu)$ times (so no U can move past
the starting point of the next U in the code).
\item Let $k(i,j)$ be the number of UR switches made using the last $j$ U's.
Then $k(i,0)=0$, for all $i$.
\item $k(i,j) \leq k(i-1,j-1)$, for all $i,j \geq 0$.
\end{itemize}
\end{theorem}
Note that this proposition implies that $c^\lambda_{\mu,\nu}$ is equal to the
number of sequences $(\mu=\mu^0, \mu^1, \mu^2, \ldots, \mu^l=\lambda)$.

\textbf{Proof:}
This theorem just follows the computational way to calculate Littlewood-Richardson
coefficients, with the only difference being that we use different notation.
The sequences $(\mu=\mu^0, \mu^1, \mu^2, \ldots, \mu^l=\lambda)$ are in 1-1
correspondence to the semistandard Young tableaux of shape $\lambda / \mu$
with 1's in the boxes in $\mu^1/\mu^0$, 2's in the boxes in $\mu^2/\mu^1$, \ldots,
and $i$'s in the boxes in $\mu^i/\mu^{i-1}$ for all $1 \leq i \leq l = u_1(\nu) = l(\nu)$.
The restriction that no U can move past the next U means that for
all $i$, $\mu^i/\mu^{i-1}$ is a horizontal $n$-strip, so the corresponding Young
tableau is indeed semistandard.  The number of $i$ boxes is the number of boxes
added to get from $\mu^{i-1}$ to $\mu^i$, which is equal to the total number of UR to RU
switches made in this step, which is $r_i(\nu)$.  This means that the Young tableau
obtained has shape $\lambda / \mu$ and weight $\nu$.  The requirement $k(i,j) \leq k(i-1,j-1)$
means that the number of $i$'s in the first $j$ rows have is less than the
number of $(i-1)$'s in the first $(j-1)$ rows for all $i$ and $j$.  This is
equivalent to saying that the reverse-row word is a lattice permutation.
\hfill $\Box$

To better illustrate the correspondence between sequences of partitions of the
form $(\mu=\mu^0, \mu^1, \mu^2, \ldots, \mu^l=\lambda)$ with the preceding conditions
and semistandard Young tableaux, we give the following example.

\begin{example}
Consider the following:  $\lambda = (4,3,2)$, $\mu = (2,1)$, and $\nu = (3,2,1)$, and
the sequence $\mu^0=\mu=(2,1)$, $\mu^1=(4,1,1)$, $\mu^2=(4,3,1),$ and $\mu^3=\lambda=(4,3,2)$.
It is straightforward though tedious to verify that this sequence does satisfies the
above conditions and hence contributes to $c^\lambda_{\mu,\nu}$.  If we follow the
algorithm in the proof of the theorem and put $i$'s in each box in $\mu^i/\mu^{i-1}$,
we get the following semistandard Young tableau:
\begin{center}
{\LARGE \young(\hfil\hfil11,\hfil22,13)}
\end{center}
We can also understand this using only the codes of these partitions.  Using the
algorithm for finding such a sequence, we would find the codes of these partitions
(not the partitions themselves) and have the following sequence (omitting leading
U's and trailing R's): $\mu^0=\mu=$ RURU, $\mu^1=$ RUURRRU, $\mu^2=$ RURRURU,
$\mu^3=\lambda=$ RRURURU.  To get from the code of $\mu^0$ to the code of $\mu^1$ the
rightmost U has to move past two R's (since the number of R's between this U and the
next rightmost U increases by two).  This means that we have to add two boxes to the
first row of $\mu$ in the first step, which are represented in the semistandard Young
tableau with 1's.  Similarly, the second U from the right does not have to move past
any R's, so there are no boxes added to the second row in the first step thus there
are no 1's in the second row in the tableau.  Again, the third U from the right must
move past one R, so one box is added in the third row and is represented by a 1 in
the third row of the tableau.  Repeating this same proceedure to get from $\mu^1$ to
$\mu^2$ gives us the boxes added in the second step which are represented by 2's in
the tableau.  Continuing in this way we can find the same semistandard Young tableau
using only the codes of the partitions.
\end{example}

Using the codes of partitions we can realize the \textit{Pieri Rules} in a new way.

\begin{corollary}\label{pieri}
(The Pieri Rules)
\begin{enumerate}
\item If $\nu = (n)$, then
$$s_\mu s_\nu = \sum_\lambda c^\lambda_{\mu,\nu} s_\lambda = \sum_{\lambda} s_{\lambda},$$
where the sum is over all $\lambda$ such that $\lambda / \mu$ is a horizontal $n$-strip.\\
\item If $\nu = (1^n) = (1,1,\ldots,1)$, then
$$s_\mu s_\nu = \sum_\lambda c^\lambda_{\mu,\nu} s_\lambda = \sum_{\lambda} s_{\lambda},$$
where the sum is over all $\lambda$ such that $\lambda / \mu$ is a vertical $n$-strip.
\end{enumerate}
\end{corollary}

\textbf{Proof:}
For part (1), $s_\mu s_\nu = \sum_{(\mu=\mu^0,\mu^1=\lambda)} s_\lambda$, where $\lambda$
is obtained from the code of $\mu$ by moving the U's to the right a total of $n$ places,
with no U moving past the starting point of the next U in the code.  Thus for any R in
the code of $\mu$, at most one U is moved past this R.  But since the number of U's moved
past the $i^\text{th}$ R from the left is the number of boxes added to the $i^\text{th}$ column,
this implies that no two of the added boxes are above each other, so $\lambda / \mu$ is
a horizontal $n$-strip.  Since $l = u_1(\nu) = 1$, each sequence has length 2, so the third
condition in the theorem, $k(i,j)\leq k(i-1,j-1)$, is satisfied trivially.
For each $\lambda$, the multiplicity of $s_\lambda$ in the summation is the number of
sequences $(\mu=\mu^0,\mu^1=\lambda)$ which is one.  Therefore $c^\lambda_{\mu,\nu} = 1$
if $\lambda$ is a horizontal $n$-strip, and $c^\lambda_{\mu,\nu} = 0$ otherwise.

For part (2), $s_\mu s_\nu = \sum_{(\mu=\mu^0,\mu^1,\ldots,\mu^n=\lambda)} s_\lambda$,
where $\mu^i$ is obtained from the code of $\mu^{i-1}$ by switching one UR to RU,
and $k(i,j)\leq k(i-1,j-1)$, for all $i,j \geq 0$.  This restriction on $k(i,j)$ implies
that the U moved to get from $\mu^{i}$ to $\mu^{i+1}$ is left of the U moved to get
from $\mu^{i-1}$ to $\mu^i$.  But since the number of R's moved past the $i^\text{th}$ U
from the right is the number of boxes added to the $i^\text{th}$ row of $\mu$,
this implies that no two of the added boxes are in the same row, so $\lambda / \mu$
is a vertical $n$-strip.  For each $\lambda$ such that $\lambda / \mu$ is a
vertical $n$-strip, the only way for the sequence $(\mu=\mu^0,\mu^1,\ldots,\mu^n=\lambda)$
to end with the partition $\lambda$ is for the rightmost box in $\lambda / \mu$ to
be added first, then the next furthest right, and so on.  Since there is only one
way to do this, the multiplicity of $s_\lambda$ in the summation is one.
Therefore $c^\lambda_{\mu,\nu} = 1$ if $\lambda$ is a vertical $n$-strip,
and $c^\lambda_{\mu,\nu} = 0$ otherwise.
\hfill $\Box$

\end{document}